\documentclass[11pt]{article}
\usepackage{latexsym,amsthm,verbatim,ifthen,graphicx,color,mathrsfs}
\usepackage[all]{xy}
\usepackage{color}

\usepackage{color}

\usepackage{lscape}
 \usepackage[ansinew]{inputenc}
 \usepackage{multicol}
 \usepackage[all]{xy}\xyoption{poly}
 \usepackage{mathrsfs}
\usepackage[psamsfonts]{eucal}
\usepackage{amssymb, amsmath}
\usepackage{geometry}
\usepackage{enumerate}
\usepackage{float}
\newtheorem{theorem}{Theorem}[section]

\newtheorem{corollary}[theorem]{Corollary}
 
 \newtheorem{proposition}[theorem]{Proposition}
 \theoremstyle{definition}
 \newtheorem{definition}[theorem]{Definition}

\def\quil{{\mathscr L}}
\newtheorem*{theorem*}{Theorem}
\def\cG{{\mathscr G}}
\def\gauge{\,\cG\,}

\def\quil{{\mathscr L}}

\newcommand{\bz}{\mathbb Z}

\def\bq{\mathbb{Q}}

\def\C{\mathbb{C}}

\def\Q{\mathbb{Q}}

\def\Z{\mathbb{Z}}

\def\im{{\rm Im\,}}
\def\hL{{\widehat{\mathbb L}}}

\def\des{{s^{-1}}}

    \newcommand{\lasu}{{\mathfrak{L}}}

      \newcommand{\derr}{\operatorname{{\rm Der}}}

 \newcommand{\lib }{\mathbb{L}}

\newcommand{\catdgl}{\operatorname{{\bf dgl}}}

\newcommand{\catcdgl}{\operatorname{{\bf cdgl}}}

\newcommand{\catss}{\operatorname{{\bf sset}}}

\newcommand{\Ho}{\operatorname{{\rm Ho}}}

    \newcommand{\Der}{\operatorname{{\rm Der}}}
    \newcommand{\id}{\operatorname{{\rm id}}}
        \newcommand{\ad}{\operatorname{{\rm ad}}}

 \newcommand{\MC}{\operatorname{{\rm MC}}}
  \newcommand{\aut}{\operatorname{{\rm aut}}}

\newcommand{\mc}{\operatorname{\MC}}
\newcommand{\wmc}{{\widetilde{\rm MC}}}

   \newcommand{\libc}{{\widehat\lib}}

\usepackage{hyperref}
\usepackage{physics}
\hypersetup{
    colorlinks=true,
    linkcolor=blue,
    filecolor=magenta,
    urlcolor=cyan,
    citecolor=blue,
}

\newcommand{\invlim}{\varprojlim}

\newcommand{\DDer}{ { \mathcal{D} \!\operatorname{er}}}
\newcommand{\dder}{\DDer}
\newcommand{\DER}{\operatorname{ \mathscr{D}\textup{er}}}

     \newcommand{\derf}{\operatorname{{\mathfrak{Der}\,}}}
     \newcommand{\der}{{\dder}}

\newtheorem{example}[theorem]{Example}				
\theoremstyle{remark}
\newtheorem{remark}[theorem]{Remark}

\begin{document}

\title{Realization of Lie algebras of derivations and moduli spaces of some rational homotopy types}

\author{Yves F\'elix, Mario Fuentes and Aniceto Murillo\footnote{The  authors have been partially supported by the MICINN grant PID2020-118753GB-I00 of the Spanish Government and the EXCEL-00827 grant of the Junta de Andaluc{\'\i}a.
 The third author also thanks the {\em Insituto de Matem\'aticas de la UNAM  en Oaxaca}, for his hospitality during a short stay in which part of this work was developed.  }}

\maketitle

\begin{abstract}
We construct Lie algebras of derivations (and identify their geometrical realization) whose Maurer-Cartan sets provide moduli spaces describing  the classes of  homotopy types of rational spaces sharing either the same homotopy Lie algebra, homology or cohomology.
\end{abstract}

\section*{Introduction}
Derivations of a Lie algebra are ubiquitous objects in topology. A particular instance is the following classical result \cite{SS1,tan}: when $L$ is a differential graded Lie algebra (dgl henceforth) characterizing  the rational homotopy type of a finite, simply connected CW-complex $X$, the dgl of positive  derivations of $L$  characterizes in the same fashion the rational homotopy type of universal covering of $B\aut^* (X)$, the classifying space of pointed self homotopy equivalences of $X$. With the recent extension of the Quillen approach to rational homotopy theory \cite{bfmt0} we were able to extend this result to  connected dgl's of derivations as long as the degree zero derivations characterizes a $\Q$-complete (in the sense of Malcev) subgroup of $\aut^*(X)$ \cite{ffm} (see also the recent approach of A. Berglund and T. Zeman \cite{ber2} to the description of the rational homotopy type of the classifying spaces of  self homotopy equivalences).

At this stage is convenient to remark that, under the mentioned extension of Quillen theory, which is the one considered in this paper, only dgl's that are complete are susceptible of being topologically realized (see \S1 for a brief compendium on this theory). Nevertheless, the reader may find other classes of dgl's whose topological realizations have been considered. See for instance, the integration procedure of the class of {\em absolute} dgl's recently developed by Roca I Lucio in \cite{ro}.

 Complete dgl's  contain much more geometrical data than their  connected cover. For instance, the Maurer-Cartan set of a dgl modulo the gauge relation ($\wmc$ set from now on) corresponds to the set of path connected component in which the realization of the given dgl decomposes.  In this paper we try to collect this extra data for some  sub Lie algebras of derivations of a given dgl,  which are complete  and still provide important geometrical information. In all these cases, their $\wmc$ sets, or the space of orbit of a certain action on them, turns out to be a moduli space governing  classes of rational homotopy types sharing certain structures.

To begin with, we consider in  \S\ref{modulis} an extended dual ``Lie version'' of the deep result of M. Schlessinger and J. Stasheff  in  \cite[Main Theorem 4.1]{SS1} which will also be considered later. Let $\pi$ be a complete connected graded Lie algebra and let $\Ho\catss_{\pi}$  be the class of  homotopy types of rational simplicial sets whose  homotopy Lie algebra is isomorphic to ${\pi}$. If $L$ denotes the  bigraded minimal Lie model of $\pi$, which is properly introduced in Theorem \ref{bigraded}, we prove the following (see Theorem \ref{tiposhomotopia} and Corollary \ref{topologico} for precise and detailed statements):

\begin{theorem*} There exists a complete sub dgl $\derf L$ of $\derr L$ such that
$$
\Ho\catss_\pi\cong\widetilde\mc(\derf L).
$$
\end{theorem*}
Via this bijection, the quotient stack $\widetilde\mc(\derf L)=\mc(\derf L)/\exp(\derf_0L)$ can be seen as the moduli space of $\Ho\catss_\pi$.

\medskip

It is important to remark that, in the simply connected case, this result was already sketched  by D. Blanc in \cite[\S3]{blanc} and explicitly developed by M. Zawodniak in his thesis \cite{za}.

\medskip

Then, in \S\ref{homo}, we construct a complete dgl of derivations which provides a moduli space governing the class $\Ho\catss^1_H$ of homotopy types of rational, finite dimensional, simply connected complexes sharing the same reduced homology with no additional structure. For it, let $\lib(V)$ be the free Lie algebra generated by $V=s^{-1}H$ and consider the dgl $L=(\lib(V),0)$ with trivial differential. With this notation, Corollary \ref{coro} can be summarized as  follows:

\begin{theorem*} There exists a complete sub dgl $\der L$ of $\derr L$ and a natural action  of  $\aut(V)$ on  $\wmc(\der L)$ for which
\begin{equation}\label{revi}
\Ho\catss^1_H\cong \wmc(\der L)/\aut(V).
\end{equation}
 Moreover,
$$
\langle \der L\rangle= \coprod_{X\in\Ho\catss^1_H}\,\,\coprod_{{\mathcal O}_X }\,\, B\aut^*_{\mathcal H}(X).
$$
\end{theorem*}
Here: $\langle\,\cdot\,\rangle$ denotes the realization functor on complete dgl's (see \S\ref{prelimi}); ${\mathcal O}_X$ denotes the (cardinality of the) orbit by the action of $\aut(V)$ of any element in $\wmc(\der L)$ representing $X$ by the bijection (\ref{revi}); and finally, $\aut^*_{\mathcal H}(X)$ is the subgroup of pointed homotopy equivalences of $X$ which induces the identity on homology.

In other words, the realization of $\der L$ is the disjoint union of simplicial sets, one for each $X\in\Ho\catss^1_H$. Moreover, each of these pieces also decomposes in as many path components as points in the orbit ${\mathcal O}_X$, and each of which is of the homotopy type of the classifying space $B\aut_{\mathcal H}^*(X)$.

Thus, in this case, $\wmc(\der L)$ is too big to describe $\Ho\catss^1_H$. Nevertheless, there is an action of $\aut(L)$ on $\mc(\der L)$ which provides the quotient stack $\mc(\der L)/\aut(L)$ responsible of $\Ho\catss^1_H$.

We remark that this result is a particular instance of the extended version in Theorem \ref{prin1}.

Finally, in \S\ref{cohomo}, we consider the augmentation ideal $A$ of a given simply connected, finite dimensional, commutative graded algebra and denote by $\Ho\catss^1_A$ the class of homotopy types of rational simply connected spaces sharing $A$ as rational (reduced) cohomology algebra. We then present a different description of $\Ho\catss^1_A$ than the one given by M. Schlessinger and J. Stasheff in \cite[Main Theorem 4.1]{SS1}. For it (see \S5 for details), denote by $L=\quil(A^\sharp)$ the classical Quillen functor on the coalgebra  given by the dual of $A$. This is a dgl with a purely quadratic differential for which  we prove (see Theorem \ref{prin2} for a precise statement):

\begin{theorem*}  There exists a complete sub dgl $\DER L$ of $\derr L$ and a natural action of $\aut(A)$ on  $\wmc(\DER L)$ so that
$$
{\Ho\catss^1_A}\cong\wmc(\DER L)/\aut(A).
$$
Moreover,
$$
\langle \DER L \rangle \simeq \coprod_{X\in {\Ho\catss^1_A}} \, \coprod_{{\mathcal O}_X} \,\, B\aut^*_{\mathcal H}(X).
$$
\end{theorem*}
Here, ${\mathcal O}_X$ denotes again the orbit by the action of $\aut(A)$ of any element in $\wmc(\DER L)$ representing $X$ by the bijection in (i). On the other hand, as before, $\aut^*_{\mathcal H}(X)$ stands for the subgroup of pointed self homotopy equivalences of $X$ which induce the identity on homology. As $X$ is rational this trivially coincides the group of self homotopy equivalences inducing the identity on cohomology.

As a consequence we can also exhibit a particular quotient stack over $\mc(\DER L)$ as a moduli space of  $\Ho\catss^1_A$.

\medskip

To prove the above results we need some technical statements which are contained in \S\ref{rea2}. This section
 extends and reformulates some results in \cite[\S6]{ffm} to obtain certain complete sub Lie algebras of a general $\derr L$ containing the whole connected cover.

\bigskip

\noindent{\bf Acknowledgement.}  We  thank the referee for his/her helpful suggestions and corrections which have considerably improved the content and presentation of this paper.

\section{Preliminaries}\label{prelimi}
This section is devoted to recall the basic facts we shall use from the homotopy theory of complete differential graded Lie algebras for which we refer to the monograph \cite{bfmt0}, or the original references \cite{bfmt3,bfmt1}, for a detailed presentation.

All considered differential graded vector spaces,  possibly endowed with additional structures, are rational and graded over $\bz$. The {\em suspension} and {\em desuspension} of  such a graded vector space $V$ is denoted by $sV$ and $\des V$ respectively. That is $(sV)_n=V_{n-1}$ and $(\des V)_n=V_{n+1}$ for any $n\in\bz$.

We often do not distinguish  objects of the  category  $\catss$  of simplicial with the topological spaces given by their realization which are therefore of the homotopy type of CW-complexes.

We denote by $\catdgl$ the category of  differential graded Lie algebras (dgl's henceforth). A dgl $L$, or $(L,d)$ if we want to specify the differential, is  {\em connected} if $L=L_{\ge 0}$.

A {\em Maurer-Cartan} element, or simply MC element, of a given dgl $L$ is an element $a\in L_{-1}$ satisfying the Maurer-Cartan equation $da=-\frac{1}{2}[a,a]$. We denote by $\mc(L)$ the set of MC elements in $L$. Given $a\in \mc(L)$, we denote by $d_a= d+\ad_a$ the {\em perturbed differential} on $L$ where $d$ is the original one and $\ad$ is the usual adjoint operator. The {\em component} of $L$ at $a$ is the connected sub dgl $L^a$ of $(L,d_a)$ given by
$$
L^a_p=\begin{cases} \ker d_a&\text{if $p=0$},\\ \,\,\,L_p&\text{if $p>0$}.\end{cases}
$$

 The {\em derivations} $\derr L$ of a given   dgl $L$ is a dgl with the usual Lie bracket and differential $D=[d,-]$:
$$
[\theta,\eta]=\theta\circ\eta-(-1)^{|\theta||\eta|}\eta\circ\theta,
\quad
D\theta=d\circ\theta-(-1)^{|\theta|}\theta\circ d.
$$

A {\em filtration} of a dgl $L$  is a decreasing sequence of differential Lie ideals,
$$L=F^1\supset\dots \supset F^n\supset F^{n+1}\supset\dots$$ such that $[F^p,F^q]\subset F^{p+q}$ for $p,q\geq 1$. In particular, the lower central series of $L$,
$$
L^1\supset\dots\supset L^n\supset L^{n+1}\supset\dots,
$$
where $L^{1}= L$ and $L^{n}= [L, L^{n-1}]$ for $n>1$, is a filtration for any dgl which satisfies $L^n \subset F^n$ for any $n\ge 1$ and  any other filtration $\{F^n\}_{n\ge 1}$ of $L$.

A {\em complete differential graded Lie algebra}, cdgl henceforth,  is a dgl $L$ equipped with a  filtration $\{F^n\}_{n\ge 1}$ for which the  natural map
$$
L\stackrel{\cong}{\longrightarrow}\varprojlim_n L/F^n
$$
is a dgl isomorphism. A cdgl {\em morphism} between cdgl's is a dgl morphism which preserves the filtrations.
We denote by $\catcdgl$ the corresponding category. By a {\em complete graded Lie algebra}, cgl hereafter, we mean a cdgl endowed with the trivial differential.

If  $L$ is a dgl  filtered by $\{F^n\}_{n\ge 1}$, its {\em completion} is the dgl
$$
\widehat L=\varprojlim_nL/F^{n}
$$
which is always complete with respect to the filtration
$$
\widehat{F}^n=\ker ( \widehat L \to L/F^n).
$$
If no specific filtration is given, the completion of a generic dgl is always taken over the lower central series. In particular, if $\lib(V)$ denotes the free Lie algebra generated by the graded vector space $V$, the completion of a dgl of the form $(\lib(V),d)$ is the cdgl
$$
\libc(V)=\varprojlim_n\lib(V)/\lib(V)^n.$$
This is an important object in this theory whose main properties are detailed in \cite[\S3.2]{bfmt0}. Note that, if $V=V_{>0}$, then $\libc(V)=\lib(V)$.

 Given a cdgl $L$  the group $L_0$, endowed  with the Baker-Campbell-Hausdorff product (BCH product henceforth), acts  on the set  $\mc(L)$ by
$$
x\,\cG\, a=e^{\ad_x}(a)-\frac{e^{\ad_x}-1}{\ad_x}(d x)= \sum_{i\geq 0} \frac{\ad_x^i(a)}{i!} - \sum_{i\geq 0} \frac{\ad_x^i(dx)}{(i+1)!},\quad x\in L_0,\,a\in\mc(L).
$$
This is the {\em gauge} action and we denote by $\widetilde\mc(L)=\mc(L)/\cG$ the corresponding orbit set. A homotopical description of the gauge action is given in  \cite[\S5.3]{bfmt0}.

The homotopy theory of cdgl's lies in the existence of a pair of adjoint functors \cite[Chapter 7]{bfmt0}, {\em (global) model} and {\em realization},
\begin{equation}\label{pair}
\xymatrix{ \catss& \catcdgl \ar@<1ex>[l]^(.50){\langle\,\cdot\,\rangle}
\ar@<1ex>[l];[]^(.50){\lasu}\\}.
\end{equation}

     The set of $0$-simplices of $\langle L\rangle$ coincides with $\mc(L)$. Moreover, if $\langle L\rangle^a$ denotes  the path component of $\langle L\rangle$ containing the MC element $a$, we have:
\begin{equation}\label{componentes}
     \langle L\rangle^a\simeq \langle L^a\rangle,\quad \langle L\rangle\simeq  \amalg_{a\in \widetilde{\mc}(L)} \,\, \langle L^a\rangle.
\end{equation}

If $L$ is connected, and for any $n\ge 1$, we have group isomorphisms
     $$
     \pi_n\langle L\rangle\cong H_{n-1}(L)
     $$
   where the group structure in $H_0(L)$ is considered with the BCH product. Under the homotopy equivalence $\langle L\rangle\simeq\mc_\bullet(L)$  between the realization of $L$ and the {\em Deligne-Getzler-Hinich groupoid} of $L$, see \cite[\S11.4]{bfmt0}, this is the original explicit isomorphism of A. Berglund $\pi_n\mc_\bullet(L)\cong H_{n-1}(L)$ in \cite[Theorem 1.1]{ber1}

     We will also use the fact  that the realization of a cdgl is invariant under perturbations. That is, for any cdgl $L$ and any $a\in \mc(L)$,
 \begin{equation}\label{ind}
     \langle L\rangle\cong\langle (L,d_a)\rangle.
\end{equation}

    Finally, the realization functor coincides with any other known geometrical realization of cdgl's. In particular, if $L$ is a $1$-connected dgl of finite type  then, see \cite[Corollary 11.17]{bfmt0}, $\langle L\rangle$ has the homotopy type of the classical Quillen realization of $L$ \cite{qui}.

     On the other hand, see again \cite[Chapter 7]{bfmt0} for details, the global model $\lasu_X$ of a simplicial set $X$ completely reflects its simplicial structure. In particular,
 the  $0$-simplices of $X$ are the Maurer-Cartan elements of $\lasu_X$.

     If $X$ is a simply connected simplicial set of finite type and $a$ is any of its vertices, then \cite[Theorem 10.2]{bfmt0}, $\lasu_X^a$ is quasi-isomorphic to $\lambda (X)$ where $\lambda$ is the classical Quillen dgl model functor \cite{qui}. Moreover, see \cite[Theorem 11.14]{bfmt0}, for any connected simplicial set $X$ of finite type, $
     \langle \lasu_X^a\rangle$ is weakly homotopy equivalent to $\bq_\infty X$ the Bousfield-Kan $\bq$-completion of $X$ \cite{BK}. Recall that, whenever $X$ is nilpotent, $\bq_\infty X$ and has the homotopy type of $X_\bq$, the rationalization  of $X$.

The category $\catcdgl$ has a {\em cofibrantly generated model structure}, see \cite[Chapter 8]{bfmt0}, for which the functors in (\ref{pair}) become a Quillen pair. With this structure the induced functors in the respective homotopy categories extend the classical Quillen equivalence between rational homotopy types of simply connected simplicial sets and homotopy types of simply connected dgl's.

A {\em model} of a connected cdgl $L$ is a connected cdgl of the form $(\libc(V),d)$ together with a quasi-isomorphism  (and hence a weak equivalence)
$$
(\libc(V),d)\stackrel{\simeq}{\longrightarrow} L.
$$
If $d$ is decomposable we say that $(\libc(V),d)$ is the {\em minimal model} of $L$ and is unique up to cdgl isomorphism.

\begin{definition}\label{minimamodel} Let $X$ be a connected simplicial set and $a$ any of its vertices. The {\em minimal model of $X$} is the minimal model of  $\lasu_X^a$.
\end{definition}

If $(\libc(V),d)$ is the minimal model of $X$ then, see \cite[Proposition 8.35]{bfmt0}, $sV\cong \widetilde H_*(X;\bq)$ and, provided $X$ of finite type, $sH_*(\libc(V),d)\cong\pi_*(\bq_\infty X)$. Again, the group $H_0(\libc(V),d)$ is  considered  with the BCH product. If $X$ is simply connected, the minimal  model of $X$ is isomorphic to its classical Quillen minimal model for which we refer to \cite{nei} or \cite{qui}.

 \section{Complete  Lie algebras of derivations}\label{rea2}

Derivations of a cdgl are essential objects in this paper. However,  even if $L$ is $1$-connected, $\Der L$ may fail to be complete and thus, their $\mc$ set are not defined and they are unable to be topologically realized as described in the past section. For instance, let $L=(\lib(x,y),0)$ with $\abs{x}=\abs{y}=2$, and consider $\theta_1,\theta_2,\theta_3\in\Der_0L$ defined by:
$$\theta_1(x)=x,\;\theta_1(y)=-y,
\;\;\; \;\;\;
\theta_2(x)=y,\; \theta_2(y)=0,\;\;\; \;\;\;
 \theta_3(x)=0,\; \theta_3(y)=x.$$
Note that
$[\theta_1,\theta_2]=-2\theta_2$, $[\theta_1,\theta_3]=2\theta_3$ and $[\theta_2,\theta_3]=-\theta_1$. Hence, for any
given filtration $\{F^n\}_{n\geq 1}$ of $\Der L$, $\theta_i\in F^n$ for any $n$ and any $i$. That is, these derivations live in the kernel of the natural map $\Der L\to \invlim_{n\geq 1} \Der L/F^n$ and thus $\Der L$ is not complete.

 Nevertheless, for any complete sub dgl $M$ of $(\derr L,D)$ we shall use the following general facts:
\begin{equation}\label{dife}
\mc(M)=\{\delta\in M_{-1}\, \,\text{such that $d+\delta$ is a differential in  $L$}\,\}.
\end{equation}
Moreover, the gauge relation is characterized by the following result:

\begin{proposition}\label{gauge} Two Maurer-Cartan elements $\delta,\eta\in\mc(M)$ are gauge related if and only if there exists an isomorphism of the form
$$
e^\theta\colon(L,d+\delta)\stackrel{\cong}{\longrightarrow}(L,d+\eta)
$$
with $\theta\in M_0$. Moreover the gauge action is given by $\theta\,\cG\,\delta=\eta$.
\end{proposition}
The proof is an obvious extension of \cite[Theorem 4.31]{bfmt0} to any complete sub dgl of derivations:
\begin{proof} Suppose first that $\delta$ and $\eta$ are gauge related. Thus, there exists $\theta\in M_0$ such that
$$
\eta=e^{\ad_\theta}(\delta)-\frac{e^{\ad_\theta}-1}{\ad_\theta}(D\theta).
$$
As $D\theta=[d,\theta]$,
$$
\frac{e^{\ad_\theta}-1}{\ad_\theta}(D\theta)=\sum_{i\geq 0} \frac{\ad_\theta^i}{(i+1)!} [d,\theta]=- \sum_{i\geq 1} \frac{\ad^i_\theta}{i!} d.
$$
Therefore,
$$
d+\eta=d+e^{\ad_\theta}(\delta)+\sum_{i\geq 1} \frac{\ad^i_\theta}{i!}( d) =e^{\ad_\theta}(d+\delta).
$$
We then use the general formula
$
e^{\ad_\theta}(d+\delta)= e^\theta(d+\delta)e^{-\theta}
$ (see for instance \cite[Proposition 4.13]{bfmt0}))
to conclude that
$$
d+\eta=e^\theta(d+\delta)e^{-\theta},\,\, \text{that is,}\,\,\,(d+\eta)e^\theta=e^\theta(d+\delta),
$$
and $e^\theta$ is the
required isomorphism.

For the other implication simply reverse the above argument.
\end{proof}
Due to this fact we often identify $M_0$ with
$$
\exp(M_0)=\{e^\theta,\,\,\theta\in M_0\},
$$
and denote  $\wmc(M)=\mc(M)/\exp(M_0)$.

If $M$ is of finite type  choose basis  $\{ \partial_i\}_{i=1}^s$  and
 $\{\sigma_\ell\}_{\ell=1}^r$
of $M_{-1}$ and $M_{-2}$ respectively and write
$$[\partial_i,\partial_j]=\sum_{\ell} \lambda_{ij}^{\ell}\, \sigma_\ell,\quad  \lambda_{ij}^{\ell}\in \Q.$$
Then, given $\delta\in M_{-1}$, the derivation $d+\delta=\sum_{i} \alpha_i\partial_i$  is a differential if and only if
$$\sum_{i,j} \lambda_{ij}^{\ell}\, \alpha_i\alpha_j=0,\quad \ell=1,\dots, r.$$
In other words, if we denote by $\mathbf{V}_L\subset\C^s$ the affine algebraic variety defined by the polynomials $\sum_{i,j} \lambda_{ij}^{\ell}\,\alpha_i\alpha_j$, with $\ell=1,\dots, r$, we conclude that
\begin{equation}\label{variedad}
\mc(M)=\{\text{rational points of $\mathbf{V}_L$}\}.
\end{equation}
so that $\wmc(M)=\mc(M)/\exp(M_0)$ can be considered as a quotient stack.

Next, consider $L=(\libc(V),d)$ a connected minimal cdgl in which $V$ is bounded above, that is $V_{>m}=0$ for some $m$. We then identify some complete sub dgl's of  $\derr L$ which conserve its ``connected cover''.   For it,  choose an arbitrary finite filtration of  $V$ by graded vector subspaces:
\begin{equation}\label{filtra}V=V^0\supset V^1 \supset\dots \supset V^{q-1}\supset V^q=0.
\end{equation}

As in \cite[\S6]{ffm},
for $n\geq 1$ and $p\geq 0$, write
$$
\hL^{n,p}(V)=\text{Span} \{ \bigl[v_1,[v_2,[\dots,[v_{n-1},v_n]\bigr]\dots\bigr]\in \hL^{n}(V), v_i\in V^{\alpha_i} \text{ and } \sum_{i=1}^n\alpha_i=p\},
$$
and define
$$
F^{n,p}=\hL^{n,p}(V)\oplus \hL^{\geq n+1}(V),$$
so that
$$
\hL(V)=F^{1,0}\supset\dots\supset F^{1,q-1}\supset F^{ 2,0} \supset\dots\supset F^{2,2q-1}\supset
\dots
\supset F^{n,0}\supset\dots\supset F^{n,nq-1}\supset
\dots
$$
In the order given by this sequence
$F^{n,p}$ takes the position
$$
t=q+\dots+(n-1)q+p+1=\frac{(n-1)nq}{2}+p+1
$$
and we define $F^t=F^{n,p}$ for $n,p$ and $t$ as above. In \cite[Proposition 6.3]{ffm} it is proved that $\{ F^t\}_{t\geq 1}$ is a filtration of $L$ for which it is complete.

This filtration of $L$ naturally determines a  decreasing sequence of sub dgl's of $\Der L$,
\begin{equation}\label{sequence}
{\mathcal F}^1\supset\dots\supset {\mathcal F}^n\supset  {\mathcal F}^{n+1}\supset \dots
\end{equation}
where, for any $n\ge1$,
$$
 {\mathcal F}^n =\{ \theta\in \Der L,\,\, \theta(F^r)\subset F^{n+r}\, \text{ for all } r\geq 0\}.$$

Note that $\{\mathcal F^n\}_{n\ge 1}$ is a filtration of the dgl $\mathcal F^1$. Moreover, a simple inspection shows that
\begin{equation}\label{gotica2}\mathcal F^1=\{\theta\in\derr L,\,\,\theta_*(V^i)\subset V^{i+1}\,\,\text{for all $i$}\},
\end{equation}
where $\theta_*\colon V\to V$ denotes the linear part of $\theta$. Then:

\begin{proposition}\label{propo}
$\mathcal F^1$ is a complete dgl.
\end{proposition}
\begin{proof}
As  $\cap_n \mathcal F ^n=0$ the map $\mathcal F^1\to \invlim_n \mathcal F^1/\mathcal F^n$ is injective.
On the other hand write a given element of $\invlim_n \mathcal F^1/\mathcal F^n$ as a  series $\sum_{n}\theta_n$ with $\theta_n\in \mathcal F ^n$. Note that, for each $v\in V$ and any integer $m\ge 1$, the series
$\sum_n \theta_n(v)$
contains a finite sum of elements in  $\lib^{m}(V)$ and thus  $\sum_n\theta_n(v)$ is a well-defined element in $\hL(V)$. Hence $\sum_{n}\theta_n\in\mathcal F^1$ and above map is also surjective.
\end{proof}
 We now ``enlarge'' the cdgl $\mathcal F^1$ as much as possible in positive degrees: starting from the original filtration (\ref{filtra}) we define a new filtration of $V$ as follows:
$$
\begin{aligned}
V&\supset V_1^1\oplus V_{\geq 2}\supset V_1^2\oplus V_{\geq 2} \supset \dots \supset V_1^{q-1}\oplus V_{\geq 2} \supset  V_{\geq 2}\supset
\\
&
\supset
 V_2^1\oplus V_{\geq 3} \supset V_2^2\oplus V_{\geq 3}\supset \dots
\supset V_2^{q-1} \oplus V_{\geq 3} \supset
V_{\geq 3}\supset \\
& \dots \\
&
\supset
V_m^1\supset V_m^2\supset \dots V_m^{q-1}\supset 0
\end{aligned}
$$
where $m$ is such that $V_{>m}=0$.
If we rename this filtration of  subspaces of $V$  by
$$
V={\mathcal V}^0\supset {\mathcal V }^1\supset{\mathcal V}^2\supset\dots\supset {\mathcal V}^{m(q-1)}\supset 0,
$$
it clearly satisfies the following property:
\begin{equation}\label{property}
 \text{${\mathcal V}^{i}_\ell\not=0\,\,$ implies  $\,\,V_{>\ell}\subset {\mathcal V}^{i+1}$.}
 \end{equation}
\begin{definition}\label{dergotica}
For this new filtration of $V$, the procedure above determines again a decreasing sequence of sub dgl's of $\derr L$ as in (\ref{sequence}), whose first term we denote by $\der L$.

\end{definition}
By Proposition \ref{propo}, $\der L$ is complete and,  in view of (\ref{gotica2}), it can be written as
$$
\der L=\{\theta\in\derr L,\,\,\theta_*({\mathcal V}^i)\subset {\mathcal V}^{i+1}\,\,\text{for all $i$}\}.
$$
Furthermore, from the characterization of $\mathcal F^1$   in (\ref{gotica2}) one easily observes that
\begin{equation}\label{gotica3}
\begin{aligned}
&\mathcal F^1_{>0}\subset \der_{>0} L,\\
&\mathcal F^1_{0}= \der_{0} L,\\
&\mathcal F^1_{<0}\supset \der_{<0} L.
\end{aligned}
\end{equation}
Moreover, we have:
\begin{proposition}\label{descripcion}
 $$\dder_k L=
 \left\{
 \begin{array}{cl}
 \Der_k L, & \text{ if } k>0,
 \\
\theta\in \Der_0 L, \textup{ such that } \theta(V^i)\subset V^{i+1}\oplus \hL^{\geq 2}(V), & \text{ if } k=0,
 \\
\theta\in \Der_k L, \textup{ such that } \theta(V)\subset \hL^{\geq 2}(V), & \text{ if } k< 0.
 \end{array}\right.
 $$
 \end{proposition}
That is, $\der L$ is a cdgl consisting of all derivations in positive degrees, those derivations of degree $0$ which increase the original filtration degree on $V$ modulo decomposables, and all derivations of negative degrees which increase the word length.

\begin{proof}

 Let $k>0$ and $\theta\in \derr_k L$. Then, for any $i$ and any non zero element of degree  $v\in{\mathcal V}^i_\ell$, it follows by (\ref{property}) that $\theta_*(v)\in V_{k+\ell}\subset {\mathcal V}^{i+1}$. By (\ref{gotica2}), $\theta\in \der_k L$.

Let  $k<0$ and $\theta\in \dder_k L$ such that such that $\theta(V)\subset \hL^{\geq 2}(V)$. By definition $\theta\in\der_k L$.   Conversely, let $\theta\in \der_k L$   and let $v\in V$ a non zero element. Assume  $v\in V_\ell$ and let $i$ be the maximal filtration index such that $v\in {\mathcal V}^i$ but $v\notin {\mathcal V}^{i+1}$. Then,   $\theta_*(v)=0$. Otherwise $ \theta_*(v)\in {\mathcal V}^{i+1}_{<\ell}$. Hence, by (\ref{property}), $V_\ell\subset {\mathcal V}^{i+2}$ which contradicts the fact that $v\notin V^{i+1}$.

Finally, for $k=0$, the obvious fact ${\mathcal F}^1_0=\der_0 L$ in (\ref{gotica3}) amounts to the required equality.
\end{proof}

\begin{remark}\label{rem} Of special interest in what follows is the particular instance of choosing the trivial filtration $V=V^0\supset V^1=0$ on $V$. In this case,
$$
\dder_k L=
 \left\{
 \begin{array}{cl}
 \Der_k L, & \text{ if } k>0,
 \\
\theta\in \Der_k L, \textup{ such that } \theta(V)\subset \hL^{\geq 2}(V), & \text{ if } k\le 0.
 \end{array}\right.
$$
\end{remark}

\section{Rational homotopy types with prescribed homotopy Lie algebras and their moduli space} \label{modulis}

 In this section we check that the method for building the moduli space of rational, simply connected homotopy types with prescribed homotopy lie algebra, already sketched in \cite[\S3]{blanc} and explicitly developed in  \cite{za}, also works in the non-simply connected case by means of the homotopy theory of cdgl's.

First, a simple inspection shows that the procedure to obtain the bigraded model of a simply connected graded Lie algebra  (see \cite[Theorems 1]{ha} or \cite[Chapter I]{ou}), dual of the classical commutative context \cite[\S3]{hasta}, extends mutatis mutandis to any connected complete cdgl:

\begin{theorem}\label{bigraded} (Complete bigraded Lie model). Let ${\pi}$ be a connected cgl. Then the cdgl $({\pi},0)$ admits a Lie minimal model
$$
\rho\colon (\libc(V),d)\stackrel{\simeq}{\longrightarrow} ({\pi},0)
$$
satisfying:

\begin{itemize}

\item $V=\oplus_{p,q\ge 0}V_p^q$ is bigraded being the lower grading the usual homological one. This bigradation extends bracket-wise to  $\libc(V)$.

    \item $dV^0=0$ and $d(V^{n+1})\subset \libc(V^{\le n})^n$, for $n\ge 0$. In particular $d$ decreases by one the upper degree so that $H(\libc(V),d)=\oplus_{p,q\ge 0}H^q_p(\libc(V,d)$ is also bigraded.

        \item $\rho\colon \libc(V^0)\twoheadrightarrow {\pi}$ is surjective, $\rho (V^n)=0$ for $n\ge 1$, $H^0(\rho)\colon H^0(\libc(V),d)\stackrel{\cong}{\longrightarrow} {\pi}$ is an isomorphism, and $H^+(\libc(V),d)=0$.

\end{itemize}

\end{theorem}

For completeness we include here the following:

\begin{proof}[Sketch of the proof]

Let ${\pi}$ be filtered by $\{F^n\}_{n\ge 1}$ so that ${\pi}\cong \varprojlim_n {\pi}/F^n$ and consider the projection $q\colon {\pi}\to {\pi}/[{\pi},{\pi}]$ onto the indecomposables of ${\pi}$. Define $V^0$ to be a space of generators of ${\pi}$ by $V^0= {\pi}/[{\pi},{\pi}]$ and choose $\rho\colon V^0\to {\pi}$ a section of $q$. Set $dV^0=0$ and extend  $\rho$ first to $\lib(V^0)\to {\pi}$ and then, by completion, to
$$\rho\colon \libc(V^0)=\varprojlim_n\lib(V^0)/\lib^n(V^0)\to \varprojlim_n {\pi}/{\pi}^n\to \varprojlim_n {\pi}/F^n={\pi}.$$
Next, define $V^1$ to be a space of relations of ${\pi}$ by $V^1= \ker\rho/[\libc(V^0),\ker\rho]$, set $\rho(V^1)=0$  and extend $d$ to $V^1$ as a section of the projection $\ker\rho\twoheadrightarrow V^1$.

For $n\ge 1$ define $V^{n+1}=H^n(\libc(V^{\le n})/[H^n(\libc(V^{\le n}),H^0(\libc(V^{\le n})]$, set $\rho(V^{n+1})=0$ and define $d\colon V^{n+1}\to \libc(V^{\le n})^n\cap \ker d$ to be a section of $\libc(V^{\le n})^n\cap \ker d\twoheadrightarrow V^{n+1}$.

\end{proof}

\begin{definition}\label{bigradeddef} The cdgl $(\libc(V),d)$ is the (complete) {\em bigraded model of ${\pi}$}. We say that the elements of $\libc(V)_p^n$ have {\em weight $p-n$}. Note that the differential $d$ preserves weight as $dV_p^n\subset \libc(V)_{p-1}^{n-1}$
\end{definition}

We now show that any cdgl whose homology is isomorphic to the cgl ${\pi}$ has a Lie model (not minimal in general) obtained by perturbing in a particular way the bigraded Lie model of ${\pi}$. The  following is again a straightforward extension to the complete connected setting of \cite[Theorems 2]{ha} or \cite[Chapter II]{ou} which is in turn the dual of \cite[Theorem 4.4]{hasta} in the commutative context.

\begin{theorem}\label{filtered}(Complete filtered Lie model).
Let $\rho\colon(\libc(V),d)\stackrel{\simeq}{\to} {\pi}$ be the bigraded model for the cgl ${\pi}$ and let $L$ be a cdgl whose homology is isomorphic to ${\pi}$.   Then, there is a Lie model of $L$ of the form,
$$
\phi\colon(\libc(V),d+\varphi)\stackrel{\simeq}{\longrightarrow} L,
$$
such that $\varphi$ increases the weight and $[\phi(v)]=\rho(v)$  for each $v\in V^0$.

Moreover, if $\gamma\colon (\libc(V),d+\psi)\stackrel{\simeq}{\longrightarrow} L$ is another Lie model under the same conditions, there exists an isomorphism
$$
f\colon  (\libc(V),d+\varphi)\stackrel{\cong}{\longrightarrow}(\libc(V),d+\psi)
$$
such that $f-\id_{\libc(V)}$ increases the weight and $\gamma f$ is homotopic to $\phi$.\hfill$\square$
\end{theorem}

\begin{definition} Let ${\pi}$ be a connected cgl and let $(\libc(V),d)$ be its bigraded model. Define the sub Lie algebra
$$\derf \libc(V) \subset\derr\,\libc(V)
$$
of derivations which raise the weight. That is, if  $W^m\subset\libc(V)$ denotes the subspace of elements of weight $m$, then $\theta\in\derf\libc(V)$ if $\theta(W^m)\subset W^{\ge m+1}$ for all $m\in\bz$.
\end{definition}

 We can now easily prove the dual of \cite[Theorem 4.1]{SS1}:

\begin{theorem}\label{tiposhomotopia}  $(\derf\libc(V),D)$ is a cdgl whose $\widetilde\mc$ set is in bijective correspondence with the set $\Ho\catcdgl_{\pi}$ of homotopy types of cdgl's whose homology is isomorphic to ${\pi}$.
\end{theorem}

\begin{proof}
Filter $\derf\libc(V)$  by $\{F^n\}_{n\ge 1}$ where,
$$
F^n=\{\theta\in \derf\libc(V),\,\, \theta(W^m)\subset W^{m+n}\,\,\text{for all $m$}\}.
$$

A simple inspection shows that $\{F^n\}_{n\ge 1}$ is indeed a filtration of the dgl $(\derf \libc(V),D)$. Moreover, $\cap_{n\ge1}F^n=0$ so that the natural map $\zeta\colon \derf\libc(V)\to\varprojlim_n \derf \libc(V)/F^n$ is injective.

On the other hand, write any $\theta\in \varprojlim_n \derf\libc(V)/F^n$ of degree $q$ as
$$
\theta=\sum_{n\ge1}\theta_n,\qquad\theta_n\in F^n,
$$
and observe that, for any $p,m\ge 0$, $\theta_n(V_p^m)=0$ as long as $n> q+m$. Hence, for any $v\in V$, $\sum_{n\ge1}\theta_n(v)$ is always a finite sum. That is, $\theta$ is a well defined element in $\derf\libc(V)$ and thus $\zeta$ is also surjective. This shows that $(\derf\libc(V),D)$ is a complete dgl. Note that $d\notin\derf\libc(V)$ as it does not rise the weight.

We next see that
  $$
  \exp\bigl(\mathfrak{Der}_0\,\libc(V)\bigr)=\{f\in\aut\libc(V)\,\,\,\text{such that $f-\id_{\libc(V)}$ raises the weight}\,\}.
  $$
Indeed, given $\theta\in\der_0\,\libc(V)$, $e^\theta-\id_{\libc(V)}=\sum_{n\ge 1}\frac{\theta^n}{n!}$ which clearly raises the weight. Conversely, given $f\in\aut\libc(V)$ such that $f-\id_{\libc(V)}$ raises the weight, the linear map
  $$
  \theta\colon V\longrightarrow \libc(V),\quad\theta(v)=\sum_{n\ge1}(-1)^{n+1}\frac{(f-\id)^n}{n},
  $$
is well defined and clearly raises the degree. In fact, the same argument used above shows that for any $p,m\ge 0$ and any $v\in V_p^m$, $(f-\id)^n(v)= 0$ for $n$ big enough. To conclude, extend $\theta$ as a derivation in $\mathfrak{Der}_0\,\libc(V)$ so that $\theta=\log f$, or equivalently, $f=e^\theta$.

  Finally, regard the $\mc$ set as in (\ref{dife}) and consider the map
  $$
 \mc\bigl(\derf\libc(V),D)\bigr){\longrightarrow}
 \Ho\catcdgl_{\pi},\qquad
  \overline\varphi\mapsto \text{homotopy type of $(\libc(V),d+\varphi)$}.
  $$
 It
   clearly factors through the orbit set,
   $$
    \wmc\bigl(\derf\libc(V),D)\bigr)=\mc\bigl(\derf\libc(V),D)\bigr)/ \exp\bigl(\derf_0\libc(V)\bigr)\stackrel{}{\longrightarrow}
 \Ho\catcdgl_{\pi},
 $$
 and, by a direct application of Theorem \ref{filtered}, this is a bijection.
\end{proof}

\begin{corollary}\label{topologico} Let ${\pi}$ be a finite type, connected cgl and let $(\libc(V),d)$ be its bigraded model. Then, the set $\Ho\catss_{\pi}$  of homotopy types of rational simplicial sets whose  homotopy Lie algebra is isomorphic to ${\pi}$ is in bijective correspondence with $\widetilde\mc\bigl(\derf\libc(V),D)\bigr)$.
\end{corollary}

\begin{proof} We first note that any  rational simplicial set  whose homotopy Lie algebra is isomorphic to ${\pi}$ is a  nilpotent, finite type simplicial set. Indeed, every complete finite type Lie algebra ${\pi}$ is degree-wise nilpotent \cite[Proposition 5.2]{ber1}. That is, for each degree $n$ there is an integer $k\ge 1$ such that any bracket of length $k$ and  degree $n$ vanishes. Moreover, if ${\pi}$ is connected, being degree-wise nilpotent is equivalent to  ${\pi}_0$ being nilpotent and acting nilpotently on ${\pi}_n$ for all $n\ge 1$. Hence, any simplicial set having ${\pi}$ as homotopy Lie algebra is necessary rational, nilpotent and of finite type.

On the other hand, the pair of adjoint functors in (\ref{pair}) restrict to equivalences between the homotopy categories of rational, nilpotent, finite type simplicial sets and that of connected cdgl's whose homology is complete and of finite type, i.e., degree-wise nilpotent. \cite[Chapter 10]{bfmt0}. To finish apply Theorem \ref{tiposhomotopia}
\end{proof}
\begin{remark}\label{moduliuno}
Identifying as in (\ref{variedad}), $\mc\bigl(\derf\libc(V),D)\bigr)$ with the rational points of the variety $\mathbf{V}_L$ with $L=(\libc(V),d)$, the above corollary exhibit the quotient stack $\mathbf{V}_L/\exp\bigl(\derf_0\libc(V)\bigr)$ as a moduli space of $\Ho\catss_{\pi}$.

 We are aware that, as we work over the rationals, this topological space  is not a quotient of a variety. Nevertheless, following \cite[\S7]{SS1}, where the authors study the commutative dual context,  one could properly define and study $\derf L$ as a scheme and $\exp\bigl(\derf_0\libc(V)\bigr)$ as an algebraic group acting on $\derf L$. In this way $\Ho\catss_\pi$ would become a quotient stack. This remark applies to the subsequent sections

\end{remark}

\section{Rational homotopy types with prescribed homology and their moduli space}\label{homo}
We describe the geometrical realization of the cdgl's of derivations provided in the past section and interpret their $\wmc$ sets from the topological point of view.

\begin{definition} Let $H$ be a simply connected graded vector space bounded above. Denote by $\Ho\catss^1_H$  the class of homotopy types of rational simply connected simplicial sets with reduced homology isomorphic to $H$. To avoid excessive notation we will not distinguish a simplicial set from the homotopy type that it represents.
\end{definition}

We fix such a graded vector space $H$ and a finite filtration of it,
$$
H=H^0\supset H^1 \supset\dots \supset H^{q-1}\supset H^q=0.
$$
This induces a filtration on  $V=s^{-1}H$ as in (\ref{filtra}). Let $L=(\lib(V),0)$ and consider the cdgl $\der L$ given in Proposition \ref{descripcion} corresponding to this filtration.

For each $X\in\catss^1_H$ denote by $G$ the subgroup of homotopy classes of self homotopy equivalences of $X$ which rise the degree of the homology filtration:
$$
G=\{[f]\in\mathcal E(X),\,\, H(f)(H^i)\subset H^{i+1}\,\,\text{for all $i$}\,\}.
$$
Consider also the subgroup $\aut^*_{G}(X)\subset \aut^*(X)$ of pointed homotopy automorphisms whose homotopy classes (free or pointed as $X$ is simply connected) live in $G$:
$$
\aut^*_{G}(X)=\{f\in\aut^*(X),\,\,[f]\in G\}.
$$
In $\catss^1_H$ there is a particular element that we denote by $X_0$ whose minimal model is $L$. This is the (co)formal space with free rational homotopy Lie algebra generated by $H$ consisting of a wedge of rational spheres, one for each generator of $H$.

\begin{theorem}\label{prin1}  (i) There are actions of $\aut(L)$ and $\aut(L)/\exp(\der_0L)$ on $\mc(\der L)$ and $\wmc(\der L)$ respectively which induce bijections
$$
\wmc(\der L)/\bigl(\aut(L)/\exp(\der_0L)\bigr)\cong\mc(\der L)/\aut(L)\cong\Ho\catss^1_H.
$$
(ii) Moreover,
$$
\langle \der L\rangle= \coprod_{X\in\Ho\catss^1_H}\,\,\coprod_{{\mathcal O}_X }\,\, B\aut^*_{G}(X).
$$
\end{theorem}
Here, ${\mathcal O}_X$ denotes the (cardinal of the) orbit by the action of $\aut(L)/\exp(\der_0L)$ of any element in $\wmc(\der L)$ providing $X$ via the bijection in {\em (i)}.
In other words, when $H$ is finite dimensional, the realization of $\der L$ is the disjoint union of simplicial sets, one for each $X\in\Ho\catss^1_H$, and each of which with as many path components as points in the orbit ${\mathcal O}_X$. Finally, each of these path components has the homotopy type of the classifying space $B\aut_{G}^*(X)$ which is nilpotent but clearly not simply connected.
\begin{remark} Note that in \emph{(ii)} we may replace the zero differential on $L$ by any other decomposable differential $d$. Indeed, in view of (\ref{dife}), any such differential  is an MC element in $\der L$ and by (\ref{ind}),
$$
\langle\der L\rangle=\langle(\der L,0)\rangle\simeq \langle(\der L,0_d)\rangle=\langle(\der L,D)\rangle
$$
where $D=[d,-]$, the  differential induced by $d$.
\end{remark}
\begin{proof}
\emph{(i)} In view of Proposition \ref{descripcion}, the $\mc$ elements of $\der L$ are simply decomposable differentials on $\lib(V)$. Therefore,
 the group $\aut(L)$ acts on $\mc(\der L)$ by
\begin{equation}\label{accion}
\varphi \cdot\delta=\varphi\delta\varphi^{-1},\quad\varphi\in\aut(L),
\quad\delta\in\mc(\der L).
\end{equation}
That is, $\varphi\cdot\delta=\delta'$ if
$$
\varphi\colon (\lib(V),\delta)\stackrel{\cong}{\longrightarrow}(\lib(V),\delta')
$$
is a dgl isomorphism. Note also that the map
$$
\mc(\der L)\to \Ho\catss^1_H,\qquad \delta\mapsto \langle (\lib(V),\delta)\rangle,
$$
induces a map on the orbit set,
\begin{equation}\label{moduli}
\mc(\der L)/\aut(L)\stackrel{\cong}{\longrightarrow}\Ho\catss^1_H,
\end{equation}
which is clearly a bijection.

On the other hand, and although $\exp(\der_0L)$ is not in general a normal subgroup of $\aut(L)$, we still can consider the short exact sequence of pointed sets:
\begin{equation}\label{fundshort}
\exp(\der_0 L)\longrightarrow\aut( L)\longrightarrow\aut(L)/\exp(\der_0L),
\end{equation}
and observe that the action of $\aut (L)$ on $\mc(\der L)$
restricts to the gauge action of $\der_0L$ on $\mc(\der L)$: $\theta\,\cG\,\delta=\delta'$ if, again,
$$
(\lib(V),\delta)\stackrel{e^\theta}{\longrightarrow}(\lib(V),\delta')
$$
is a dgl isomorphism.

 Hence, $\aut(L)/\exp(\der_0L)$ acts on the orbit set $\mc(\der L)/\exp(\der_0L)=\wmc(\der L)$ and
$$
\wmc(\der L)/\bigl(\aut(L)/\exp(\der_0L)\bigr))\cong\mc(\der L)/\aut(L).
$$
This and (\ref{moduli}) proves \emph{(i)}

\emph{(ii)} By (\ref{componentes}), the number of connected components of $\pi_0\langle \der L\rangle$ is in bijective correspondence with $\wmc(\der L)$. But, in view of
 \emph{(i)}, each homotopy type of $\Ho\catss^1_H$ contains as many $\wmc$ elements of $\der L$ as point in $\mathcal O_X$.  Hence, the number of connected components of $\langle\der L\rangle$ is as asserted .

Next, choose $d\in\wmc(\der L)$ which again corresponds to a decomposable differential $d$ in $L=\libc(V)$. Then, the (algebraic) component $(\der L)^d$ is the connected cdgl,
 $$(\der L)^d_k=
 \left\{
 \begin{array}{cl}
 \Der_k L & \text{ if } k>0,
 \\
\der_0 L\cap\ker D, & \text{ if } k=0,
 \end{array}\right.
 $$
 whose differential is $D=[d,-]$, induced by $d$.  By \cite[Theorem 7.13]{ffm}, if we denote by $X\in\catss^1_H$ the (rational homotopy type of the) simplicial set whose minimal model is $(\lib(V),d)$, we deduce that
$$
\langle (\der L)^d\rangle\simeq B\aut^*_G(X)
$$
and \emph{(ii)} follows.
\end{proof}
 The following particular instance is of special interest. If we choose in $H$ the trivial filtration $H=H^0\supset H^1=0$, Theorem \ref{prin1} reads:
 \begin{corollary}\label{coro}
(i) There are  actions of $\aut(L)$ and $\aut(V)$ on $\mc(\der L)$ and $\,\wmc(\der L)$ respectively which induce bijections
$$
\wmc(\der L)/\aut(V)\cong\mc(\der L)/\aut(L)\cong\Ho\catss^1_H.
$$

\bigskip

(ii) Moreover,
$$
\langle \der L\rangle= \coprod_{X\in\Ho\catss^1_H}\,\,\coprod_{\mathcal O_X }\,\, B\aut^*_{\mathcal H}(X).
$$
\end{corollary}
Here, for each $X\in \Ho\catss^1_H$, $\mathcal H$ denotes the subgroup  of homotopy classes of self homotopy equivalences that induce the identity on homology. Again, ${\mathcal O}_X$ denotes the (cardinality of the) orbit by the action of $\aut(V)$ of any element in $\wmc(\der L)$ representing $X$ by the bijection in {\em (i)}.

\begin{proof}
Recall from Remark \ref{rem} that in this case
$$
\dder_0 L=\{\theta\in \derr_0 L, \textup{ such that } \theta(V)\subset \hL^{\geq 2}(V)\}.
$$
 Then,
 $$
 \exp(\der_0L)=\{\varphi\in\aut(L),\,\,\varphi_*=\id_V\},
 $$
 where $\varphi_*\colon V\to V$ denotes the induced map on the indecomposables. Hence, (\ref{fundshort}) becomes
\begin{equation}\label{bien}
 \exp(\der_0 L)\longrightarrow\aut( L)\longrightarrow\aut(V))
\end{equation}
 and Theorem \ref{prin1}\emph{(i)} translates to \emph{(i)}. With this, and the fact that $G=\mathcal H$ in this case,  \emph{(ii)} is obvious.
\end{proof}

\bigskip

In view of the isomorphism in (\ref{moduli}),
$$\Ho\catss^1_H\cong\MC(\dder L)/\aut(L),$$
 we can identify the quotient stack $\MC(\dder L)/\aut(L)$ as a moduli space of the set of simply connected rational homotopy types with prescribed reduced homology $H$.
Moreover, two proportional (non trivial) differentials in $\mc(\der L)$ are in the same orbit. Hence, as the polynomials defining $\mathbf{V}_L$ are homogeneous, we can think of $\Ho\catss^1_H-\{X_0\}$ as a quotient stack of a subvariety of a projective space.

 \begin{example} Let $H$ be the vector space with two generators of degrees 2 and 4 and another two generators of degree 6. We compute the moduli space of $\Ho\catss^1_H$.

 For it,
let $L=(\lib(V),0)$ where $V$ is the vector space with generators $x,y,z,w$ of degrees 1, 3, 5 and 5 respectively.   We endow $V$ with the trivial filtration.
 Then, $\der_{-1}L$ is a 3-dimensional vector space, generated by the derivations $\delta_y$, $\delta_z$ and $\delta_w$ defined by
 $$\delta_y(y)=[x,x],\quad \delta_z(z)= [x,y],\;\;\; \delta_w(w)=[x,y]$$
 and are zero otherwise. In this particular case, $\mc(\der L)=\der_{-1}L$. Moreover, one easily checks that the gauge action is trivial and thus, in view of Corollary \ref{coro}\emph{(i)},
  $$
  \mc(\der L)/\aut(L)\cong \mc(\der L)/\aut(V)\cong\Ho\catss^1_H.
  $$
Hence,
If we use $\{\delta_y,\delta_z,\delta_w\}$ as basis,  we identify 4 different orbits in $\mc(\der L)/\aut(V)$ represented by the derivations $(0,0,0)$, $(\alpha,0,0)$ with $\alpha\neq 0$, $(0,\beta,\gamma)$ with either $\beta$ or $\gamma$ not zero, and $(\alpha,\beta,\gamma)$ with $\alpha\neq 0$ and either $\beta$ or $\gamma$ not zero. By considering to which spaces correspond these differentials we obtain:
 $$
 {\Ho\catss^1_H}=\{\,S^2\vee S^4\vee S^6\vee S^6,\quad (S^2\times S^4)\vee S^6,\quad \C P^2\vee S^6\vee S^6,\quad \C P^3\vee S^6\,\}
 $$
In other words the moduli space of $\Ho\catss^1_H$ consists of 4 points, where $S^2\vee S^4\vee S^6\vee S^6$ is a closed point, $\C P^2 \vee S^6$ is an open point, and any neighborhood of $(S^2 \times S^4)\vee S^6$ or $\C P^2\vee S^6\vee S^6$ contains the point $\C P^3\vee S^6$. As a finite topological space, it is characterized by its corresponding poset in which $x\le y$ if and only if $x$ belongs to the closure of $y$.

In the following figure we depict this poset and the algebraic variety $\mathbf{V}_L$ which in this case is all $\C^3$, or $\C P^2$ if we consider the corresponding projective variety by removing the origin. There, we identify the rational points belonging to each orbit of the moduli space: the origin is the only point in its orbit and corresponds to $S^2\vee S^4\vee S^6\vee S^6$; all rational points of the $\delta_y$ axis are in the orbit of $\C P^2\vee S^6\vee S^6$; rational points of the plane generated by $\{\delta_z,\delta_w\}$ conform the orbit of $(S^2\times S^4)\vee S^6$; and the rest of rational points are in the orbit of $\C P^3\vee S^6$.
\hskip1cm
\begin{figure}[H]
\centering
  \includegraphics[scale=0.53]{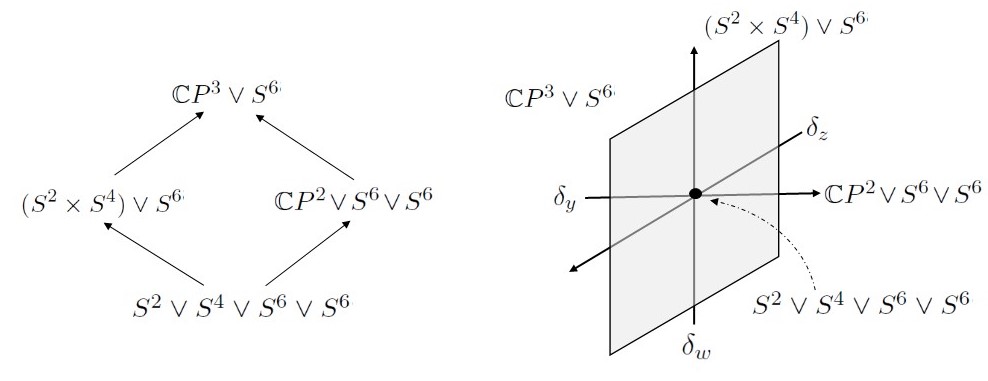}
\end{figure}
 \end{example}

Note that in general, for any $L$, the zero differential is always the only point in its orbit space $\mc(\der L)/\aut (L)$, it is a closed set and it corresponds to $X_0$, a wedge of spheres determined by a set of generators of $H$.

Remark also that ${\Ho\catss^1_H}\cong\mc(\der L)/\aut (L)$ is not always a finite space as one can check by, for instance, by computing the example in which $H=\operatorname{Span}\{x,y,z\}$ with $\abs{x}=5,\abs{y}=6$ and $\abs{z}=22$.

An interesting property of the set of elliptic homotopy types sharing the same homology is:

\begin{proposition} The set $\,\mathbf{Ell}_H\subset{\Ho\catss^1_H}$  of homotopy types of elliptic spaces  is always an open subset of the moduli space.
 \end{proposition}
 \begin{proof}
 Fix $\delta\in \MC(\dder_{-1} L)$ such that $(L,\delta)$ is elliptic.  Since $L$ is of finite type the function
 $$\dim_k\colon \mc(\der L)\to \Z,\,\, \dim_k(\delta')= \dim H_k(L,\delta')
 $$
 is well defined  for each $k\ge 1$. Note that, if $V=V_{\le N}$ then $\dim_k(\delta)=0$ for $k\ge 2N$ (see for instance \cite[Corollary 1, \S32]{FHT}). Moreover, by elementary linear algebra, regarding $\dim_k(\delta')$ as $\dim\ker\delta'|_{ L_k}-\dim\im\delta'|_{L_{k+1}}$ the map $\dim_k$ is semicontinuous for all $k$.
In particular, for each $k\ge1$ there is a neighborhood $\theta_k$ of $\delta$ such that $\dim_k(\delta')\le \dim(H_k(L,\delta))$ for any $\delta'\in\theta_k$.

Consider the open set $\theta=\cap_{i=2N}^{3N}\theta_i$ in which $\dim_k=0$ for all $k=2N,\dots,3N$.  This implies that for any $\delta'\in\theta$,  $(L,\delta')$ is elliptic. Otherwise, by \cite[Theorem 33.3]{FHT}, if $(L,\delta')$ is hyperbolic,  there must be an integer $k_0$ with $2N\le k_0\le 3N$ such that $\dim_{k_0}(\delta')\not=0$.

Finally, if we denote by $p\colon\mc(\der L)\to\mc(\der L)/\aut(L)$ the projection, then $p(\theta)\subset \mathbf{Ell}_H$ is clearly an open set of the moduli space $\mc(\der L)/\aut(L)$ containing the orbit of $\delta$.

 \end{proof}

\section{Rational homotopy types with prescribed cohomology algebra and their moduli space}\label{cohomo}

Let $L=(\libc(V),d)$ be a connected free cdgl  and consider in $V$, which is supposed to be bounded above, the trivial filtration so that Remark \ref{rem} applies.

\begin{definition} Define ${\mathscr D}\text{er}(L)$ as the complete sub dgl of $\dder L$ given by:
  $$\mathscr{D}\text{er}_k L=
 \left\{
 \begin{array}{cl}
 \der_k L, & \text{ if } k\ge 0,
 \\
 \{\eta\in \Der_k L, \textup{ such that } \eta(V)\subset \libc^{\geq 3}(V) \}, & \text{ if } k< 0.
 \end{array}\right.
 $$
 \end{definition}
This cdgl will be essential in what follows.
\begin{definition}
Consider a simply connected,   commutative graded algebra of finite dimension whose augmentation ideal we denote by $A$. Define
$
\Ho\catss^1_A$ as the class  of homotopy types of rational simply connected simplicial sets
with reduced cohomology algebra isomorphic to $A$. Again, we will not distinguish a simplicial set from the homotopy type that it represents.
\end{definition}

Recall that, given $X\in \Ho\catss^1_A$,
a classical fact, see for instance \cite[III.3.(9)]{tan}, states that the differential $d$ in $\quil(A^\sharp)$, necessarily quadratic, is naturally identified with the cup product of $X$. Here $\quil$ denotes the classical Quillen functor from coalgebras to Lie algebras. We then fix $A$, rename $L=\quil(A^\sharp)$ and prove:

\begin{theorem}\label{prin2}
(i) There is an action of $aut(A)$ on $\wmc(\Der L)$ which induces a bijection
$$\wmc(\DER L)/\aut(A)\cong  {\Ho\catss^1_A}.
$$
(ii) Moreover,
$$
\langle \DER L \rangle \simeq \coprod_{X\in {\Ho\catss^1_A}} \, \coprod_{\mathcal O_X}\,\, B\aut^*_{\mathcal H}(X).
$$
\end{theorem}

Once again, ${\mathcal O}_X$ denotes the (cardinality of the) orbit by the action of $\aut(A)$ of any element in $\wmc(\der L)$ representing $X$ by the bijection in (i).

\begin{proof} Write $L=(\lib(V),d)$ where $V=s^{{-1}}A^\sharp$ and $d$ is quadratic. Recall that the differential in $\DER L$ is  $[d,-]$. Hence,
 an MC element of $\DER L$ is, by definition, a derivation $\delta$ of $L$ such that $\delta(V)\subset \lib^{\ge 3}(V)$ and $d+\delta$ is a differential on $L$. In what follows we use the following trivial identification,
$$
\mc(\DER L)\cong\{
\,
 \text{differentials $d+\delta$ on $L$ such that  $\delta(V)\subset \lib^{\ge 3}(V)$}\,\},
 $$
so that
 $$
 \mc(\DER L)\subset \mc(\der L,0)=\{\,\text{decomposable differentials on $L$}\,\}.
 $$
 We then  consider the stabilizer $\aut_d(L)$ of $\mc(\DER L)$ of the action (\ref{accion}) of the (non differential automorphisms) $\aut(L)$ on $\mc(\der L,0)$. That is:
$$\aut_{d}(L)=\{ \varphi\in \aut(L), \text{ such that the quadratic part of } \varphi^{-1}\circ d\circ \varphi \text{ is } d\}.$$

On the other hand, the surjective map
 $$\MC(\DER L ) \longrightarrow {\Ho\catss^1_A}, \quad d+\delta\mapsto \langle (L,d+\delta)\rangle,
 $$
  clearly induces a map on the set of orbits
 $$
 \MC(\DER L)/\aut_d(L) \longrightarrow {\Ho\catss^1_A}.
 $$
Now, $\langle (L,d+\delta)\rangle\simeq \langle (L,d+\delta')\rangle$ if and only if there is dgl isomorphism $$
\varphi\colon (L,d+\delta)\stackrel{\cong}{\longrightarrow}(L,d+\delta').
$$
Thus, $\varphi\in\aut_d(L)$ and $\varphi\cdot(d+\delta)=d+\delta'$. This shows that the above map is also injective and we have a bijection
$$\mc(\DER L )/\aut_d(L)\cong {\Ho\catss^1_A}.$$

Next, observe that $\exp(\DER_0 L)\subset \aut_d(L)$ and the quotient $\aut_d(L)/\exp(\DER_0 L)$ is trivially identified to the group of automorphism of $V$ which respects the quadratic differential $d$. That is,
$$
\aut_d(L)/\exp(\DER_0 L)\cong \{ \phi\colon V\stackrel{\cong}{\longrightarrow} V, \text{ such that } \phi^{-1}\circ d\circ\phi=d \}.
$$
But this group is in bijective correspondence with the algebra  automorphisms $\aut (A)$  and we have the following
 short exact sequence, analogue of (\ref{bien}),
$$
\exp(\DER_0L)\to \aut_{d}(L)\to \aut (A)
$$
Next, observe that the action of $\aut_d(L)$ on $\mc(\DER L)$ restricts to the gauge action of $\DER_0 L$ on $\mc(\DER L)$. Hence,  as in the proof of Theorem \ref{prin1}\emph{(i)}, we deduce that $\aut(A)$ acts on $\wmc(\DER L)$ and
\begin{equation}\label{final}
\wmc(\DER L)/\aut(A)\cong \mc(\DER L)/\aut_d(L)\cong {\Ho\catss^1_A}.
\end{equation}

On the other hand, via this bijection,  each homotopy type $X$ of ${\Ho\catss^1_H}$ contains as many $\wmc$ elements of $\DER L$ as points in the orbit $\mathcal O_X$
 and thus, the number of path components of $\langle \DER L\rangle$ is as asserted in \emph{(ii)} for a general $A$.

Finally, since $A$ is finite dimensional and
 $\DER_{\ge0}L=\der_{\ge 0}L$, every connected component of $\langle \DER L\rangle$ is necessary of the homotopy type of $B\aut^*_{\mathcal H}(X)$ for the corresponding $X\in{\Ho\catss^1_A}$, just as in Corollary \ref{coro}\emph{(ii)}.

\end{proof}

\begin{remark} We can also exhibit the set of simply connected homotopy types sharing the same cohomology algebra $A$ as a quotient stack. Indeed, in view of (\ref{final}),
$$
\Ho\catss^1_A\cong \mc(\DER L)/\aut_d(L)
 $$ which by (\ref{variedad}) is a quotient of rational points in a variety. Moreover, observe that
$$\dder_{-1}L=\dder^2_{-1}L\oplus \DER_{-1} L,$$
where $\dder^2_{-1}L=\{ \eta\in \Der_{-1}L, \text{ such that } \eta(V)\subset \mathbb L^2(V)\}$ are the quadratic derivations.
Therefore,  we can identify the  $\MC(\DER L)$,  with the intersection of the algebraic variety $\MC(\dder L,0)$ with the affine linear subspace $d+\DER_{-1}L$:
$$\MC(\DER L)= \MC(\dder L,0)\cap (d+\DER_{-1} L),$$
\end{remark}

 \begin{example} Consider the commutative graded algebra $A$ generated by  the elements $a,b,c,p,q$ of degrees 4, 6, 13, 15, 19 respectively, and whose only non trivial products are:
 $$a p=q= bc .
 $$
We determine the moduli space of $\Ho\catss^1_A$. Note that $L=\quil(A)=(\lib(V),d)$ where $V$ is generated by elements $x,y,z,u,v$ of degrees 3, 5, 12, 14, 18 respectively. The differential is given by
 $$dv=[x,{u}]+[y,z]$$
 and zero on any other generator.

 We now compute
 $
 \wmc(\DER L)/\aut(A)$. First, we check that $\DER_{-1}L$ is generated by 3 derivations $\delta_z,\delta_{u}$ and $\delta_{v}$defined by
 $$\delta_z(z)=\ad_x^2(y),\;\;\; \delta_{u}(u)= ad_y^2(x),\;\;\; \delta_{v}(v)=\ad_x^4(y),$$
 and zero otherwise. A direct computation shows that a general element $\alpha\delta_z+\beta\delta_{u}+\gamma\delta_{v}$ is in $\mc(\DER L,D)$ if and only if $\alpha=\beta$.

To compute the gauge action  we first check that  $\DER_0L$ is generated by three derivations defined by
 $$\theta(u)= \ad_x^3(y), \;\;\; \theta'(v)= \ad^2_x(y),\;\;\; \theta''(v)= -\ad_y^3(x),$$
 and zero otherwise. Another straightforward computation shows that
 $$D\theta=-\delta_{v},\;\;\; D\theta'=0,\;\;\; D\theta''=0\quad\text{and}\quad [\DER_0L,\DER_{-1}L]=0.
 $$
Therefore, the only non-trivial gauge action is
 $$(t\theta)\gauge (\alpha\delta_z+\alpha\delta_{u}+\gamma \delta_{v})=\alpha\delta_z+\alpha\delta_{u}+(\gamma+t)\delta_{v}$$
 for any $t\in \Q$. Therefore, in $\widetilde\MC(\DER L)$, we can take representatives with $\gamma=0$, so that
 $$\widetilde\MC(\DER L)=\{ \alpha(\delta_z+\delta_{u}), \text{ with } \alpha\in \Q\}.$$

 Finally, any automorphism  $\phi\in\aut(A)$ is given by
 $$\phi(a)=\lambda_a a,\;\;\; \phi(b)=\lambda_b b,\;\;\;
 \phi(c)=\lambda_c c,\;\;\; \phi(p)=\lambda_p p, \;\;\;
 \phi(q)=\lambda_qq$$
 where the scalars are non zero and satisfy
 $$\lambda_a\lambda_p=\lambda_q=\lambda_b\lambda_c.$$
 For $\alpha\neq 0$ choose $\phi$ with $\lambda_a=\lambda_c=\lambda_q=1/\alpha$ and $\lambda_b=\lambda_p=1$. Then, one checks that the action of $\phi$ on the $\mc$ element $\delta_z+\delta_{u}$ gives
 $\alpha(\delta_z+\delta_{u})$.

 Therefore, in $\wmc(\DER L)/\aut(A)$ there are only two orbits corresponding to $\alpha\not=0$ and $\alpha=0$. By Theorem \ref{prin2}\emph{(i)}, We conclude that
 $$
 \Ho\catss^1_A=\{X_0,X_1\}
 $$
 where $X_0= \langle (L,d)\rangle$ is the formal space of cohomology algebra $A$ and $ X_1=\langle (L,d+\delta_z+\delta_{u})\rangle $ is the rationalization of
$SU(6)/SU(3)\times SU(3)$.

Moreover, as a moduli space, $\Ho\catss^1_A$ has the Sierpinski topology in which $X_1$ is open.

 \end{example}

\bigskip
\bigskip\bigskip
\noindent {\sc Institut de Math\'ematiques et Physique, Universit\'e Catholique de Louvain, Chemin du Cyclotron 2,
1348 Louvain-la-Neuve,
         Belgique}.

\noindent\texttt{yves.felix@uclouvain.be}

\bigskip

\noindent{\sc Departamento de \'Algebra, Geometr\'{\i}a y Topolog\'{\i}a, Universidad de M\'alaga, 29080 M\'alaga, Spain.}

\noindent
\texttt{m\_fuentes@uma.es, aniceto@uma.es}

 \end{document}